\input amstex
\input epsf.tex
\documentstyle{amsppt}

\refstyle {A}
\magnification=\magstep1 

\hoffset .25 true in
\voffset .2 true in

\hsize=6.1 true in
\vsize=8.5 true in

\def\op{\operatorname}

\def\a{\alpha}
\def\b{\beta}
\def\l{\lambda}

\def\e{\epsilon}
\def\m{\mu}
\def\s{\sigma}

\def\g{\frak g}
\def\gq{{\frak g}_q}
\def\hg{\tilde\frak g}
\def\hgq{\tilde{\frak g}_q}

\def\CC{\Bbb C}

\def\Cq{\Bbb C (q)}

\def\R{\frak R}
\def\Rq{{\frak R}^q}
\def\U{\bold U}

\def\ssl{\frak {sl}}

\def\ot{\otimes}

\topmatter
\title
Solutions of the  Yang-Baxter equation\\
and quantum $\frak {sl}(2)$
\endtitle

\author Maxim Vybornov \endauthor
\address Department of Mathematics, Yale University, 
10 Hillhouse Ave,  New Haven, CT 06520 
\endaddress
\email  mv\@math.yale.edu
\endemail

\abstract We construct a quantum deformation of a family of 
the Yang-Baxter equation solutions naturally arising from
a Lie algebra $\ssl(2,\CC)$.
\endabstract
\endtopmatter

\document

\head Introduction \endhead

Compared to the well established theory of quantum groups,
the theory of quantum Lie algebras is much less developed.
Quantum Lie algebras and related topics of bicovariant
differential calculus are studied in [Be1], [Be2], 
[CSWW], [DF], [DH], [LS], and  [Ma1], 
to mention just a few papers. However, many problems still seem 
to be open.

It is known that a solution $\R$
of the Quantum Yang-Baxter Equation (QYBE) can be naturally associated
to any Lie algebra. Moreover, in some sense the Yang-Baxter relation
is equivalent to the fact that the Lie bracket satisfies Jacobi identity
and the invertibility of $\R$ is equivalent to the fact that
the Lie bracket is antisymmetric. A non-trivial quantum 
deformation $\Rq$ of $\R$ can be regarded therefore as a version of
quantum Lie algebra with the Yang-Baxter relation serving as 
quantum Jacobi identity and the invertibility serving
as quantum antisymmetry. 

In this paper we introduce a construction of a quantum deformation 
of the classical (family of) QYBE solutions in the case of 
$\ssl(2,\CC)$. The construction is presented in the framework of the
Penrose-Kauffman graphical calculus.
This example motivates an abstract definition of a 
(generalized) quantum Lie algebra
and its module. The category of such modules is shown to have a 
natural tensor structure.

The paper is organized as follows. 
In 1.1 we discuss two families of QYBE solutions
naturally arising from a Lie algebra (cf. [Ma2]). In 1.2
we discuss a family of QYBE solutions
naturally arising from a module over a Lie algebra.
In 1.3 we recall some basic facts about $\ssl(2,\CC)$
and the classical version of the Penrose-Kauffman graphical calculus,
while in 2.1 we recall similar facts about $U_q(\ssl(2,\CC))$
and the quantum graphical calculus.
In 2.2. we introduce two families of QYBE solutions 
deforming the classical families for $\ssl(2,\CC)$ and in 
2.3. we discuss quantum $\ssl_q(2,\CC)$-modules.  

The proofs of all the lemmas are straightforward calculations
or graphical calculus exercises; they are left to the 
reader.

\head 1. Classical families of the QYBE solutions \endhead

\subhead 1.1. Lie algebras \endsubhead

Let $\CC$ be the field of complex numbers, and let $\g$ be a 
$\CC$-vector space. Suppose that $\g$ is an algebra, i.e. 
there is a multiplication map
$[\ ,\ ]:\g\ot\g\to\g$, and  let $\hg=\g\oplus\CC$. We define two 
families of linear maps 
$$
\aligned 
\R(\l),\R(\l)'&: \hg\ot\hg\to\hg\ot\hg \\
\R(\l)&: (x+\a)\ot (y+\b)\mapsto (y+\b)\ot (x+\a)+ [x,y]\ot\l \\ 
\R(\l)'&: (x+\a)\ot (y+\b)\mapsto (y+\b)\ot (x+\a)+ \l\ot [x,y] 
\endaligned
$$
where $x,y\in\g$, $\a,\b,\l\in\CC$. We set 
$\R_{12}(\l)=\R (\l)\ot 1_{\hg}$ and $\R_{23}(\l)=1_{\hg}\ot\R (\l)$.

\proclaim{Lemma 1.1.1} a) $\R(\l)\R(\l)'=\R(\l)'\R(\l)=1$ for any $\l$,
if and only if $[\ ,\ ]$ is antisymmetric;

b) $\R_{12}(\l)\R_{23}(\l)\R_{12}(\l)
=\R_{23}(\l)\R_{12}(\l)\R_{23}(\l)$ for any $\l$,
if and only if $[\ ,\ ]$ satisfies Jacobi identity;

c) $\R_{12}(\l)'\R_{23}(\l)'\R_{12}(\l)'
=\R_{23}(\l)'\R_{12}(\l)'\R_{23}(\l)'$ for any $\l$,
if and only if $[\ ,\ ]$ satisfies Jacobi identity.
\endproclaim

For any two linear spaces $V$ and $W$ let $\s (V\ot W)$
be the transposition of the tensor factors 
$\s (V\ot W):V\ot W \to W\ot V$. 
Let $\R_{21}(\l)=\s(\hg,\hg)\R(\l)\s(\hg,\hg)$.

\proclaim{Lemma 1.1.2} For any $\l$,
$\R(\l)\R_{21}(-\l)=\R_{21}(\l)\R(-\l)=1$.
\endproclaim

Therefore, a classical Lie algebra $\g$ provides a one-parameter 
family of the QYBE solutions in $\op{End}\hg\ot\hg$. 

\subhead 1.2. Modules over Lie algebras \endsubhead

Let $\g$ be a Lie algebra. Let $V$ be a $\CC$-vector space.
We consider a linear map $A: \g\ot V\to V$, and define a
family of intertwiners 
$$
\aligned 
\R_V(\l)& : \hg\ot V\to V\ot\hg \\
\R_V(\l)& : (x+\a)\ot v\mapsto v\ot (x+\a)+ A(x,v)\ot\l \\
\endaligned
$$
where $x\in\g$, $v\in V$, $\a,\l\in\CC$.

\proclaim{Lemma 1.2.1} $(\R_V(\l)\ot 1_{\hg})\circ 
(1_{\hg}\ot\R_V(\l))\circ 
(\R(\l)\ot 1_V)=
(1_V\ot \R(\l))\circ (\R_V(\l)\ot 1_{\hg})\circ (1_{\hg}\ot \R_V(\l))
$ for any $\l$, if and only if A is a $\g$-action.
\endproclaim

\proclaim{Lemma 1.2.2} Let $V, W$ be two $\g$-modules. Then
for any $\l$, 
$(\s (V\ot W)\ot 1_{\hg})\circ (1_V\ot\R_W(\l))\circ 
(\R_V(\l)\ot 1_W)=(1_W\ot \R_V(\l))\circ (\R_W(\l)\ot 1_V)\circ 
(1_{\hg}\ot \s (V\ot W))$.
\endproclaim

\subhead 1.3. $\ssl(2,\CC)$ and Penrose-Kauffman graphical calculus
\endsubhead

Let us consider a special case $\g=\ssl(2,\CC)$. Let $V_n$ denote
the $(n+1)$-dimensional irreducible $\ssl(2,\CC)$-module. In this
notation the fundamental representation is $V_1$ and the 
adjoint representation is $V_2$. An {\sl intertwiner} is a morphism
in the category of $\ssl(2,\CC)$-modules.
The Lie bracket is the unique (up to a constant) 
intertwiner $[\ ,\ ]: V_2\ot V_2\to V_2$.

Penrose-Kauffman graphical calculus is a way to represent the
$\ssl(2,\CC)$-intertwiners by planar diagrams. 
We will now recall the very basic conventions of graphical calculus;
for details we refer the reader to [CFM], [KL], [FK]
and references thereof (see, however, the warning below). 
The module $V_1$ (or rather $\op{Id}_{V_1}$) is depicted by a 
solid vertical 
strand $\vcenter{\epsfbox{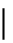}}$. 
An intertwiner $a:(V_1)^{\ot m}\to(V_1)^{\ot n}$ is depicted
by certain curves connecting $m$ distinct points on one horizontal 
line and $n$ distinct points on another horizontal line lying below
the first one. Only simple intersections are allowed. 
The module $V_n$ (or rather the Jones-Wenzl projector $p_n$) 
is depicted by a box marked $n$ 
with $n$ strands attached to its top and $n$ strands attached
to its bottom. 
$$
\vcenter{\epsfbox{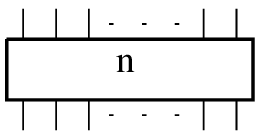}}
$$
In particular, $V_2$ is represented by the following diagram
$$
\vcenter{\epsfbox{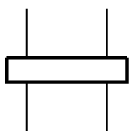}}
$$
We omit the marking of the box if it does not lead to a confusion.
A graphical version of the relation $V_1\ot V_1\simeq V_2\oplus V_0$
is depicted as follows:
$$
\vcenter{\epsfbox{box.eps}}{}={}
\vcenter{\epsfbox{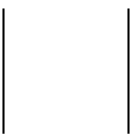}}{}+{}\frac{1}{2}{}
\vcenter{\epsfbox{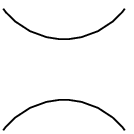}}\tag 1.3.1
$$
Using the graphical calculus we can present the Lie bracket $[\ ,\ ]$, 
the transposition $\s (V_1\ot V_1)$, and our families 
$\R(\l),\R(\l)',\R_{V_n}(\l)$ 
as follows:
$$
\matrix
[\ ,\ ]=\vcenter{\epsfbox{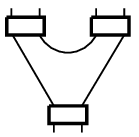}} & 
\s(V_1\ot V_1)=\vcenter{\epsfbox{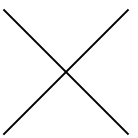}} \\
\endmatrix 
$$
$$
\matrix 
\R(\l)=\vcenter{\epsfbox{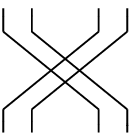}}+\l
\vcenter{\epsfbox{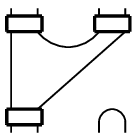}} \\
 \\
 \\
\R(\l)'=\vcenter{\epsfbox{cross.eps}}+\l
\vcenter{\epsfbox{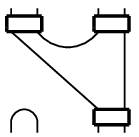}} \\
 \\
 \\
\R_{V_n}(\l)=\vcenter{\epsfbox{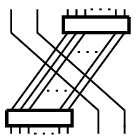}}
+\l\vcenter{\epsfbox{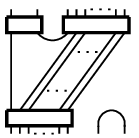}}
\endmatrix
$$

\head 2. Families of the QYBE solutions arising from quantum 
 $\ssl(2,\Bbb C)$ \endhead

\subhead 2.1. $U_q(\ssl(2,\Bbb C))$ and Penrose-Kauffman graphical 
calculus
\endsubhead

The quantum group $\U=U_q(\ssl(2,\Bbb C))$ is a certain Hopf algebra
over the field $\Cq$, defined by generators and relations.
For details on the definition of $\U$ we refer the reader 
to any of the numerous books on quantum groups (cf. [Dr], [J1], [J2]). 
We follow the notation and conventions of [FK]. 
For $n\in\Bbb Z$ there is a $\U$-module $V^q_n$ deforming the
usual $(n+1)$-dimensional module over $\ssl(2,\Bbb C)$.
A linear map  between two $\U$-modules is called an 
{\sl intertwiner} if it is a morphism in the category of $\U$-modules.
We are especially interested in two $\U$-modules: the
deformation of the natural representation $V^q_1$ and the
deformation of the adjoint representation $V^q_2$.
There are certain basic intertwining operators which we would like to
mention
$$
\matrix\format\l \\
\e_1: V_1^q\ot V_1^q\to V^q_0 \\
\e_1(v^1\ot v^1)=\e_1(v^{-1}\ot v^{-1})=0 \\
\e_1(v^{-1}\ot v^1)=1, \e_1(v^1\ot v^{-1})=-q \\
\delta_1:V^q_0\to V_1^q\ot V_1^q \\
\delta_1=v^1\ot v^{-1}-q^{-1}v^{-1}\ot v^1 \\
\check R_{11}: V_1^q\ot V_1^q\to V_1^q\ot V_1^q \\
\check R_{11}=\s R
\endmatrix
$$
where $\s$ is the transposition of the tensor factors, and 
$R(V_1^q\ot V_1^q)$ is the representation in
$V_1^q\ot V_1^q$ of a certain distinguished element $R\in \U\ot\U$,
called the universal $R$-matrix. We have the following relation
$$
\matrix
\check R_{11}=q^{\frac{1}{2}}\delta_1\circ\e_1+q^{-\frac{1}{2}}I
\endmatrix\tag 2.1.1
$$
\remark{Remark} For the explanation of the appearance of 
the square root $q^{\frac{1}{2}}$ see e.g. [FK, 1.3].
\endremark
Clearly, there is a unique (up to a constant) intertwiner
$[\ ,\ ]_q:V_2^q\ot V_2^q\to V_2^q$ which provides a structure
of an algebra for the $\Cq$-vector space $V_2^q$. We call this
intertwiner the {\sl quantum Lie bracket}. We also denote 
$V_2^q$ by $\gq$ or $\ssl_q(2,\CC)$,
and $V_1^q\ot V_1^q\simeq  V_2^q\oplus V_0^q\simeq \gq\oplus\Cq$
by $\hgq$.

The quantum version of Penrose-Kauffman graphical calculus is a way to 
represent the $\U$-intertwiners by planar diagrams. 
Again, we recall the very basic conventions of graphical calculus;
for details we refer the reader 
to [KL] or [FK] (see, however, the warning below). 
The deformation of the natural 
representation $V_1^q$ is depicted by a solid vertical strand 
$\vcenter{\epsfbox{v.eps}}$.
An intertwiner $a:(V_1^q)^{\ot m}\to(V_1^q)^{\ot n}$ is depicted
by certain curves connecting $m$ distinct points on one horizontal 
line and $n$ distinct points on another horizontal line lying below
the first one. Only simple intersections are allowed. At each 
intersection we specify the type of intersection: overcrossing or 
undercrossing. Thus, any diagram can be viewed as a projection 
of a system of curves in three dimensions.

The quantum adjoint representation  $V_2^q$ is depicted by a box
with two strands attached to its top and two strands attached to its 
bottom. We need the quantum version of the relation (1.3.1):
$$
\vcenter{\epsfbox{box.eps}}{}={}
\vcenter{\epsfbox{id.eps}}{}+{}\frac{1}{[2]}{}
\vcenter{\epsfbox{de.eps}}
$$
where $[2]=q+q^{-1}$ (in general, $[n]=\frac{q^n-q^{-n}}{q-q^{-1}}$), 
and the graphical version of the relation (2.1.1):
$$
\vcenter{\epsfbox{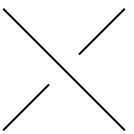}}{}=
q^{\frac{1}{2}}{}\vcenter{\epsfbox{de.eps}}+
q^{-\frac{1}{2}}{}\vcenter{\epsfbox{id.eps}}
$$
The quantum Lie bracket $[\ ,\ ]_q$ is depicted as follows
$$
\vcenter{\epsfbox{bra.eps}}
$$
Roughly speaking, to recover the classical graphical calculus from its
quantum version one has to stop distinguishing between overcrossings
and undercrossings and set $q=1$.
\remark {Warning} In our conventions the intertwiners ``go downward'' 
while in [FK] the intertwiners ``go upward''. Therefore,
all our pictures are upside down compared to those of [FK].
\endremark

\subhead 2.2. $\ssl_q(2,\CC)$ in the framework of graphical calculus
\endsubhead

We consider the following eight intertwiners 
$r,r',a,a',b,b',c,c': \hgq\ot\hgq\to\hgq\ot\hgq$:

$$\matrix 
r=\vcenter{\epsfbox{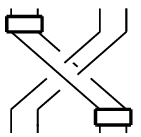}} &
a=\vcenter{\epsfbox{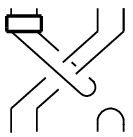}} &
b=\vcenter{\epsfbox{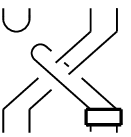}} &
c=\vcenter{\epsfbox{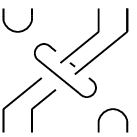}} \\
\endmatrix
$$
and
$$
\matrix 
r'=\vcenter{\epsfbox{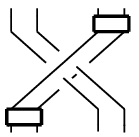}} &
a'=\vcenter{\epsfbox{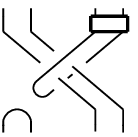}} &
b'=\vcenter{\epsfbox{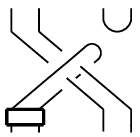}} &
c'=\vcenter{\epsfbox{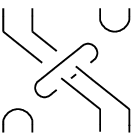}} \\
\endmatrix
$$
It is clear how to define the above intertwiners algebraically,
using the projectors and injectors (see [FK, 1.4]) and 
the basic intertwiners $\e_1,\delta_1, \check R_{11}, 
\check R_{11}^{-1}$. We leave this exercise to the reader.  
We define two families of intertwiners
$$
\aligned
\Rq(\l),\Rq(\l)'&: \hgq\ot\hgq\to\hgq\ot\hgq \\
\Rq(\l)&=r+(\l-\frac{1}{[2]})a+\frac{1}{[2]([2]\l-1)}b+\frac{1}{[2]^2}c \\
\Rq(\l)'&=r'+(\l-\frac{1}{[2]})a'+\frac{1}{[2]([2]\l-1)}b'+
\frac{1}{[2]^2}c'
\endaligned
$$
where $\l\in\Cq$, $\l\neq\frac{1}{[2]}$. It is easy to see that
$$
\matrix
\Rq(\l)& = & \vcenter{\epsfbox{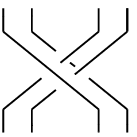}}+
\l\vcenter{\epsfbox{sa.eps}}+
\frac{\l}{[2]\l-1}\vcenter{\epsfbox{sb.eps}}\\
\endmatrix
$$
and
$$
\matrix
\Rq(\l)'& = & \vcenter{\epsfbox{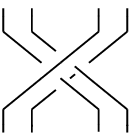}}+
\l\vcenter{\epsfbox{sla.eps}}+
\frac{\l}{[2]\l-1}\vcenter{\epsfbox{slb.eps}}\\
\endmatrix
$$

\proclaim{Lemma 2.2.1} 
a) For any $\l\neq\frac{1}{[2]}$, 
$\Rq(\l)\Rq(\l)'=\Rq(\l)'\Rq(\l)=1$;

b) For any $\l\neq\frac{1}{[2]}$, 
$\Rq_{12}(\l)\Rq_{23}(\l)\Rq_{12}(\l)
=\Rq_{23}(\l)\Rq_{12}(\l)\Rq_{23}(\l)$;

c) For any $\l\neq\frac{1}{[2]}$, 
$\Rq_{12}(\l)'\Rq_{23}(\l)'\Rq_{12}(\l)'
=\Rq_{23}(\l)'\Rq_{12}(\l)'\Rq_{23}(\l)'$.
\endproclaim

Let us specialize the family $\Rq(\l)$ to the case when 
$\l=\frac{\m}{[2](1-q^{-2})}$, $\m\in\CC$, $\m\neq 0$.
Then it is easy to see that if $q=1$, the family $\Rq(\l)$
degenerates to our classical family $\R(\m)$, $\m\neq 0$.
A similar statement is true for $\Rq(\l)'$ and $\R(\m)'$. Thus, we can
regard the $\Cq$-vector space $\hgq$ equipped with $\Rq(\l)$
as a quantum Lie algebra. The relation of Lemma 2.2.1.a) 
should be regarded as quantum antisymmetry and the relation
of Lemma 2.2.1.b) or Lemma 2.2.1.c) as quantum Jacobi identity.

\remark{Remarks} 1) We would obtain two more families of QYBE solutions 
if we turn all the diagrams in the definition of $\Rq(\l),\Rq(\l)'$ 
upside down;

2) Along the lines of Majid [Ma2], 
a solution of the QYBE can be constructed in the framework of [LS].
It is easy to see that
such a  solution arising in the case of $\ssl(2,\CC)$ is 
a member of our family.
A similar technique can be applied in the case of $\ssl(n,\CC)$;
we do not know, however, how to extend it to the case of an arbitrary 
semisimple Lie algebra.
\endremark

\subhead 2.3. $\ssl_q(2,\CC)$-modules in the framework of graphical 
calculus \endsubhead

Let $V$ be an irreducible finite dimensional $\U$-module, i.e. 
$V=V^q_n$ for some $n$. 
Let us consider the following four intertwiners 
$\hgq\ot V\to V\ot\hgq$:
$$
\matrix
r(V)=\vcenter{\epsfbox{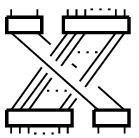}} &
a(V)=\vcenter{\epsfbox{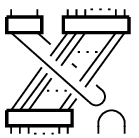}} &
b(V)=\vcenter{\epsfbox{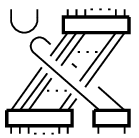}} &
c(V)=\vcenter{\epsfbox{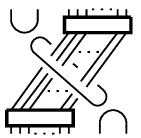}} \\
\endmatrix
$$
It is clear how to define the above intertwiners algebraically. 
We leave this exercise to the reader.  
We define a family of intertwiners
$$
\aligned
\Rq_V(\l)&: \hgq\ot V\to V\ot\hgq \\
\Rq_V(\l)&= r(V)+(\l-\frac{1}{[2]})a(V)+\frac{1}{[2]([2]\l-1)}b(V)+
\frac{1}{[2]^2}c(V)
\endaligned
$$
where $\l\in\Cq$, $\l\neq\frac{1}{[2]}$. It is easy to see that
$$
\matrix
\Rq_V(\l)& =  \vcenter{\epsfbox{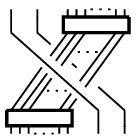}}+
\l\vcenter{\epsfbox{qan.eps}}+
\frac{\l}{[2]\l-1}\vcenter{\epsfbox{qbn.eps}}\\
\endmatrix
$$
\proclaim{Lemma 2.3.1} For any $\l\neq\frac{1}{[2]}$,
$(\Rq_V(\l)\ot 1_{\hg})\circ (1_{\hg}\ot\Rq_V(\l))\circ 
(\Rq(\l)\ot 1_V)=
(1_V\ot \Rq(\l))\circ (\Rq_V(\l)\ot 1_{\hg})\circ (1_{\hg}\ot \Rq_V(\l))
$.
\endproclaim

For two $\U$-modules $V$, $W$ we define $R(V\ot W)$ as the 
representation of the universal $R$-matrix in $V\ot W$.
We set $\check R(V\ot W)=\s R(V\ot W)$, where $\s$ is the 
transposition of the tensor factors.
The following Lemma is a $q$-deformation of Lemma 1.2.2. 

\proclaim{Lemma 2.3.2} Let $V, W$ be two irreducible finite
dimensional $\U$-modules. Then for any $\l\neq\frac{1}{[2]}$,
$(\check R(V\ot W)\ot 1_{\hg})\circ (1_V\ot\Rq_W(\l))\circ 
(\Rq_V(\l)\ot 1_W)=
(1_W\ot \Rq_V(\l))\circ (\Rq_W(\l)\ot 1_V)\circ (1_{\hg}\ot 
\check R(V\ot W))$. 
\endproclaim

Let us consider the $\U$-module $V^q_n$.
Let us specialize the family $\Rq_{V^q_n}(\l)$ to the case when 
$\l=\frac{\m}{[n](1-q^{-2})}$, $\m\in\CC$, $\m\neq 0$.
Then it is easy to see that if $q=1$, the family $\Rq_{V^q_n}(\l)$
degenerates to our classical family $\R_{V_n}(\m)$, $\m\neq 0$.

\remark{Remarks} 1) We would obtain another family of intertwiners 
if we turn all the diagrams in the definition of $\Rq_V(\l)$ 
upside down.

2) Let us consider the case when $V=V_2^q\gq$. We can regard 
$a(\gq)$ as a map $\gq\ot\gq\to\gq$ and $c(\gq)$ as a map
$\gq\to\gq$. Then $a(\gq)=[2](1-q^{-2})[\ ,\ ]_q$, 
$c(\gq)=-(q^3+q^{-3})\op{Id}_{\gq}$. Let us consider the components
of $\Rq_{12}(\l)\Rq_{23}(\l)\Rq_{12}(\l)$ and 
$\Rq_{23}(\l)\Rq_{12}(\l)\Rq_{23}(\l)$ mapping $\gq\ot\gq\ot\gq$
to $\gq\ot\Cq\ot\Cq$. These components may be regarded as
maps from $\gq\ot\gq\ot\gq$ to $\gq$. Then the Yang-Baxter equation 
implies
$$
(q^2+q^{-2}-1)[\ ,\ ]_q\circ([\ ,\ ]_q\ot 1)=
[\ ,\ ]_q\circ(1\ot [\ ,\ ]_q)\circ
(1^{\ot 3}-r(\gq)\ot 1)\tag 2.3.1
$$
It is clear that the maps $T^h$ and $\s^h$ introduced in [V]
may be considered as $\Cq$-linear maps between $\Cq$-linear spaces.
Then $T^h$ is equal to $[\ ,\ ]_q$ up to multiplication by a constant
and $1^{\ot 2}-r(\gq)=\frac{q^2+q^{-2}}{2}(1^{\ot 2}-\s^h)$.
Therefore, (2.3.1) is equivalent to the second of the relations
(6.1) of [V]. The first relation (6.1) can be obtained in a similar 
way.
\endremark
\medskip

Motivated by the above examples we give the following abstract 
definitions.

\definition{Definition 2.3.3} A $\Cq$-vector space $\gq$ is called 
a (generalized) 
quantum Lie algebra if it is equipped with a $\Cq$-linear family 
$\R(\l): \hgq\ot\hgq\to\hgq\ot\hgq$, where $\hgq=\gq\oplus\Cq$, 
of the invertible solutions of the Yang-Baxter equation.
A $\Cq$-vector space $V$ is called a (generalized)
$\gq$-module if there is a $\Cq$-linear family 
$\R_V:\hgq\ot V\to V\ot\hgq$ satisfying the relation of Lemma 2.3.1.
\enddefinition
For example, any $V_n^q$ is an $\ssl_q(2,\CC)$-module. 

\proclaim{Lemma 2.3.4} Let $V$, $W$ be two $\gq$-modules.
Then the vector space $V\ot W$ has a natural structure of 
a $\gq$-module: $\R_{V\ot W}=(1_V\ot\R_W)(\R_V\ot 1_W)$.
\endproclaim

Such a tensor product is clearly associative. Therefore, the
category of $\gq$-modules is a monoidal category.

\head  Acknowledgments \endhead 

I am very grateful to my advisor Igor Frenkel for the formulation 
of the problem and very useful discussions as well as
continuous support; and to J. Ding, 
P. Etingof, M. Khovanov,
S. Majid, A. Malkin, and N. Reshetikhin 
for numerous useful discussions.

\enddocument
\end